\documentclass [11pt]{article}

\title {The Word Problem for the Singular Braid Monoid}
\author {Oliver T. Dasbach
       \thanks{The work of the first author was supported by the
Deutsche Forschungsgemeinschaft (DFG)}
       \thanks{e-mail: \sl kasten@math.columbia.de,
       \sl http://www.math.uni-duesseldorf.de/home/kasten}\\
          {\small \sl Columbia University}\\
          {\small \sl Department of Mathematics}\\
          {\small \sl New York, NY 10027} \and
         Bernd Gemein
        \thanks{e-mail: \sl gemein@rz.uni-duesseldorf.de}\\
          {\small \sl Heinrich-Heine-Universit\"at}\\
          {\small \sl Mathematisches Institut}\\
          {\small \sl Universit\"atsstr. 1, D-40225 D\"usseldorf}
          }
\date{Beta release, V. 0.99}

\input{epsf}
\newtheorem{satz}{Satz}[section]
\newtheorem{example}[satz]{Example}

\newtheorem{lemma}[satz]{Lemma}
\newtheorem{definition}[satz]{Definition}
\newtheorem{proposition}[satz]{Proposition}
\newtheorem{theorem}[satz]{Theorem}
\newtheorem{corollary}[satz]{Corollary}

\newenvironment
  {proof}{\noindent {\bf Proof }}{$\Box$\bigskip \medskip}

\newenvironment
  {remark}{\bigskip \noindent {\bf Remark }}{\bigskip}

\newcommand {\Z} {\mathbb{Z}}

\newcommand {\pr} [2]
   { \langle #1  \vert \, #2 \rangle }

\usepackage{latexsym}
\usepackage{amssymb}

\newcommand {\s} {{\sigma}}
\renewcommand {\t} {{\tau}}
\newcommand {\cd}{\mbox{cd}} 

\begin{document}

\maketitle

\begin{abstract}
We give a solution to the word problem for the singular braid monoid 
$SB_n$. The complexity of the algorithm is quadratic in the product of the
word length and the number of the singular generators in the 
word. 
Furthermore we algebraically reprove a result of Fenn, Keyman and 
Rourke that the monoid embeds into a group and we compute the cohomological
dimension of this group.
\end{abstract}

\section{Introduction}
Back in the 20th Emil Artin introduced the braid group
\cite{Artin}.
He gave a presentation and already showed how to solve the word problem for 
this group.

When the theory of Vassiliev knot invariants started in the early 90ths, 
it became also interesting - both from the point of view of mathematics
as of physics - to look at singular braids (see e.g. \cite{Birman2,
Baez, Hutchings, FRZ, FKR, Vershinin}), where transversal 
self-intersections are allowed. These singular braids form a monoid. 

While for the word problem in the braid group many different 
solutions are known (\cite{Garside, Thurston, BKL, FGRRW}), 
for the singular braid monoid such an algorithm was not 
known and it seems very difficult to extend one of the
solutions of the word problem in the braid group  to 
the singular braid monoid.

The aim of this paper is to give an algorithm solving the word problem in the
singular braid monoid and - as Artin did for the braid group - to give
informations on the algebraical structure of the singular braid monoid. 
This will be done by using traditional algebraic
tools, such as properties of {\it HNN-}extensions of groups.

We proceed as follows. As proved by Fenn, Keyman and Rourke \cite{FKR}
the singular braid monoid $SB_n$ embeds into a group $SG_n$.
Since the proof given there involved some geometrical arguments which do 
not seem to have generalizations for much more general settings, for 
example for other ``singular Artin groups'', in the course of this text 
we will give a group theoretical proof of it.

This embedding theorem allows us to use the tools of classical group theory
for the solution of the word problem.  
We will work out the structure of a certain subgroup of finite index in 
$SG_n$ as an iterated {\it HNN-}extension - with some nice properties - of a
subgroup of the braid group. 
Britton's lemma together with the known solution to the word problem for
the braid group now allows us to give a solution to the word problem for the
singular braid monoid.

Along the path of our proof we can give some information about the
group $SG_n$. For example we will compute its cohomological
dimension to be $n-1+\lfloor n/2 \rfloor$. 
For the easiest case $SG_3$ we will construct a $K(SG_3,1)-$space and compute
the homology of the group. 

In Section \ref{section proofs} we will give some technical but
necessary proofs. The trustful reader can skip this section.

The first author would like to thank Joan Birman, Yair Glasner and 
Fritz Grunewald for many helpful discussions.

The second author thanks Wilhelm Singhof for useful suggestions and remarks.

\section{The singular braid monoid}

\begin{figure}[h] 
\begin{center}
\parbox{3.5cm}{\epsfbox{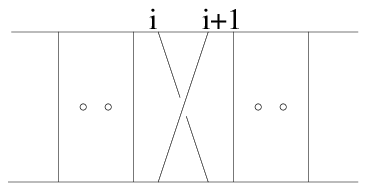}}
\parbox{3.5cm}{\epsfbox{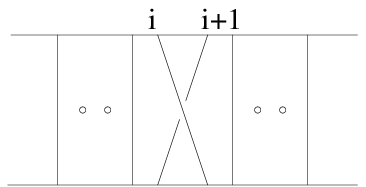}}
\parbox{3.5cm}{\epsfbox{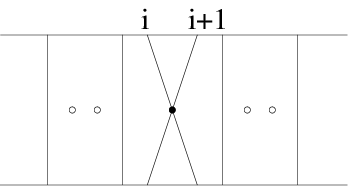}}
\end{center}
\caption{\label{elementary braids} $\sigma_i, \sigma_i^{-1}$ and $\tau_i$}
\end{figure}

The theory of Vassiliev invariants made it interesting to investigate knotted
objects having a finite number of transversal self-intersections. As such a
generalization of the braid group $B_n$ we get the singular braid 
monoid $SB_n$ generated by the elementary singular braids 
$\s_1, \dots, \s_{n-1}$ and $\t_1, \dots \t_{n-1}$
depicted in Figure \ref{elementary braids}.

Joan Birman \cite{Birman2} and independently John Baez \cite {Baez} gave a presentation
for the singular braid monoid. The generators $\t_j, \, j>1$, 
can be expressed in terms of  
$\s_1, \dots, \s_{n-1}$ and $\t_1$,  and one can show that the monoid presentation 
is equivalent to the following presentation:

\begin{proposition} [\cite{DG1}] \label{second presentation of SB}
The monoid $SB_n$ is generated by the elements $$\s_i^{\pm 1}, \, i=1, \dots, n-1,
\, \mbox{ and }
\t_1$$ satisfying the following relations: 
\begin{eqnarray}
\s_i \s_i^{-1}&= &1 \mbox{ for all } \, i \\
\s_i \s_{i+1} \s_i = \s_{i+1} \s_i \s_{i+1}
&\mbox{and}& \s_i \s_j = \s_j \s_i \quad \mbox{ for } \, \, j>i+ 1
\\
\s_2 \s_1^2 \s_2 \t_1 &=& \t_1 \s_2 \s_1^2 \s_2\\
\s_i \t_1 &=& \t_1 \s_i \quad \mbox{ for } \, i \neq 2 \\
\s_2 \s_3 \s_1 \s_2 \t_1 \s_2 \s_3 \s_1 \s_2 \t_1 &=&
\t_1 \s_2 \s_3 \s_1 \s_2 \t_1 \s_2 \s_3 \s_1 \s_2  \quad  \mbox{ for } n >  3.
\end{eqnarray}   
\end{proposition}

It would be very arduous to work with a semigroup so
in the course of this article we will heavily make use of the following embedding
theorem:

\begin{theorem}[\cite{FKR}] \label{embeds in a group}
The singular braid monoid $SB_n$ embeds into a group that will be denoted by $SG_n$.
That means $SG_n$ is the group that we get by regarding the presentation of
$SB_n$ as a group presentation.
\end {theorem}

Since the original proof of this theorem involved geometrical considerations 
and since the theorem lies on our road, in Section \ref{section embedding} we will 
give an algebraic proof of it.

For later use we will need the following theorem:

\begin{theorem}[\cite{FRZ}] \label{commutation} 
For a braid $\beta \in B_n$ the following are equivalent:
\begin{enumerate}
\item $\s_j \beta = \beta \s_k$
\item $\s_j^r \beta = \beta \s_k^r$ for some nonzero integer $r$ 
\item $\t_j \beta = \beta \t_k$
\item $\t_j^r \beta = \beta \t_k^r$ for some positive integer $r$.
\end{enumerate}
\end{theorem}

Actually,  we only need a lemma that is proved in \cite{FRZ} as an application.
Let $\eta$ be the homomorphisms 
$$
\eta: SB_n \longrightarrow \Z B_n
$$
of the singular braid monoid into the integral group ring of the braid group,
that is induced by the map: $\t_i \mapsto \s_i - \s_i^{-1}$ and $\s_i \mapsto \s_i$.

We have:

\begin{lemma} [\cite{FRZ}] \label {injective for two}
The homomorphism $\eta$ is injective for the subset $SB^{(2)}_n$ of
$SB_n$ of all singular braids having at most two singularities.
\end{lemma}

\begin{remark}
The more general conjecture of Joan Birman \cite{Birman2}, that the homomorphism 
$\eta$ is injective is still open. For a discussion of this problem see 
\cite{FRZ}, \cite{Jarai}, \cite{Zhu} and  \cite{DG1}.
\end{remark}

\section{{\it HNN-}extensions of groups}

Our main tool is the concept of {\it HNN-}extensions of groups (see e.g. 
\cite{LS} or \cite{MKS}).
Let $H = \pr {S\,}{\mbox{rel. }H \,}$ be a group with a set of generators $S$ and
relations $\mbox{rel. }H$ 
and $U$ and $V$ two isomorphic subgroups of $H$ together with
an isomorphism $\Phi$. 

The {\it HNN-}extension $G$ of $H$ relative to $U$ and $V$ is

$$
G \cong \pr {H, t \,} {\mbox{rel. }H,\, t u t^{-1} = \Phi(u), \, u \in U \,}.
$$  

The element $t$ is called {\it stable letter} and $H$ is the {\it base group}.
In our cases $\Phi$ is always the identity, so from now on we only consider such
{\it HNN-}extensions. 
By the classical result of Higman, Neumann and Neumann the group $H$ is embedded into
$G$, that means the subgroup of $G$ that is generated by the elements of $H$ is
isomorphic to $H$.

Central for our solution to the word problem for the singular braid monoid 
is the following beautiful result of Britton, often quoted as Britton's lemma:

\begin{lemma}[\cite{Britton}] \label {Brittons lemma}
Let $H= \pr {S \,} {\mbox{rel. }H \,}$ be a presentation of the group $H$ with a set of generators
$S$ and relations $\mbox{rel. }H$ in these generators.

Furthermore let $G$ be an {\it HNN-}extension of $H$ of the following form:
$$ G = \pr {S,t \, } {\mbox{rel. H}, t^{-1} u t = u, \, u \in U \, }$$ 
for some subgroup $U \subset H$.  

Let $w$ be a word in the generators of $G$ which involves $t$.
If $w=1$ in $G$ then $w$ contains a subword $t^{-1} u t$ or $t u t^{-1}$
where $u$ is a word in $S$, and $u$, regarded as an element of the group $H$,
belongs to the subgroup $U$.
\end {lemma}

We will apply Britton's Lemma to a subgroup of finite index in the singular braid 
group $SG_n$ to show that this subgroup and hence $SG_n$ itself has a solvable 
word problem.

In general, if the base group $H$ has a solvable word problem it does not mean
that an {\it HNN-}extension $G$ as in Britton's Lemma has to have one. In addition
there must be a test whether a given element of $H$ is in $U$ or not,
or equivalently whether it commutes with $t$ or not.

In our proof $t$ is always a singular generator and when $u$ lies in the braid group
$B_n$ the test whether $u$ commutes with $t$ is established by the solution to the word
problem in $B_n$ and Lemma \ref{injective for two}. 

\section{The subgroups $SGD_{n,i}$ of $SG_n$}

In general a good reference for all used facts about braid groups is \cite{Birman3}.
Our notation in the following is a modification of the notation in
\cite{Chow}. Especially we think of $B_{n-1}$ as the subgroup of
$B_n$ generated by $\{ \s_2, \dots, \s_{n-1}\}$ rather than of the one
generated by $\{ \s_1, \dots, \s_{n-2} \}$.

\begin{definition}
Let $SGD_{n,i}$ be the preimage of 
$$\Sigma_{n-i}=Sym(\{i+1,\dots,n\}) 
\subset \Sigma_n=Sym(\{1,\dots, n\})$$ 
of the natural homomorphism
\begin{equation} \label{map onto Sn}
SG_n \longrightarrow \Sigma_n
\end{equation}
 and let $D_{n,i}$ (resp. $SD_{n,i}$) be the corresponding
subgroup of $B_n$ (resp. the submonoid of $SB_n$).
Especially we have $SGD_{n,0}=SG_n$.
We will call the kernel $SGP_n$ of the homomorphism in (\ref{map onto Sn})
the pure singular braid group. So $SGP_n := SGD_{n,n-1}$.
\end{definition}

\begin {lemma}
The underlying geometry gives us an embedding 
$$
\phi_{n,i}: SG_{n-i} \longrightarrow SGD_{n,i}
$$
induced by the map
$$
\s_j \mapsto \s_{j+i}  \qquad \qquad \t_j \mapsto \t_{j+i}.
$$

The same holds for the embedding $SB_{n-i}$ into $SD_{n,i}$ and 
$B_{n-i}$ into $D_{n,i}$.
\end{lemma}

\begin {lemma}  \label{coset enumeration}
A system of Schreier right cosets of $SGD_{n,i-1}$ modulo $SGD_{n,i}$ is given
by 
$$M_{i,j}:= \s_i \s_{i+1} \cdots \s_j, \quad j=i, \dots, n-1$$
and the identity.

We get as generators for $SGD_{n,i}$:
\begin{eqnarray*}
a_{k,j}&:=& \s_k \s_{k+1} \cdots \s_j^2 \s_{j-1}^{-1}\cdots \s_k^{-1}, \qquad 1 \leq k \leq j \leq n-1, \quad k \leq i\\
X_{k,j}&:=& \s_j \s_{j-1} \cdots \s_k \t_k \s_{k+1}^{-1} \cdots \s_{j-1}^{-1} \s_j^{-1}, \qquad i \geq j \geq k \geq 1\\
\t_{i+1} &\mbox{ and }& \s_{i+1}, \dots, \s_{n-1}.\\
\end{eqnarray*}
\end{lemma}

\subsection{A presentation for the subgroup $SGD_{n,1}$}

Essentially for our considerations will be to work out a group presentation for
the subgroup $SGD_{n,1}$ of $SG_n$.
A presentation for the corresponding subgroup of $B_n$ was given by Chow:

\begin{theorem}[Chow, \cite{Chow}] \label{Chow}
The subgroup $D_{n,1}$ of $B_n$ is generated by the elements $a_{1,1},
\dots,
a_{1, n-1}$ and $\s_2, \dots, \s_{n-1}$
subject to the relations
\begin{enumerate}
\item The relations of $B_{n-1}$ generated by $\s_2, \dots, \s_{n-1}$ hold
\item
\begin {eqnarray}
\s_i a_{1,k} \s_i^{-1} &=& a_{1,k}   \quad \mbox{ for } k \neq i,
i-1\label{rel9}\\
\s_i a_{1,i} \s_i^{-1} &=& a_{1,i-1}\label{rel10}\\
\s_i a_{1, i-1} \s_i^{-1} &=& a_{1, i-1}^{-1} a_{1,i} a_{1,i-1}.
\label{rel11}
\end{eqnarray}
\end{enumerate}

Furthermore the subgroup of $D_{n,1}$ generated by $a_{1,1}, \dots,
a_{1,n-1}$ is a free subgroup of rank $n-1$ and lies normal in $D_{n,1}$.
\end{theorem}

\begin{remark}

This presentation may be simplified according to the philosophy of our
paper. With the setting $a := a_{1,1}=\sigma_1^2$ one can get:

\begin{corollary} \label{Artin group}
The subgroup $D_{n,1}$ of $B_n$ is generated by the elements $a$ and 
$\sigma_2, \dots \sigma_{n-1}$ subject to the usual braid relations 
and the relations:
\begin{eqnarray}
\s_i a &=& a \s_i \quad \mbox{ for } \, i \geq 3, \label{bcommutator} \\
\s_2 a \s_2 a &=& a \s_2 a \s_2. \label{fourbraid}
\end{eqnarray}
\end{corollary}

Corollary \ref{Artin group} shows that $D_{n,1}$ again is an Artin group.
It recently gained some new interest when tom Dieck \cite{tomDieck2} 
studied representations of it. Because of its own geometrical meaning it is 
also called cylinder braid group $ZB_{n-1}$.
\end{remark}

\begin{theorem} \label{Pres SGD_n1}
$SGD_{n,1}$ is generated by the elements $\s_2, \dots, \s_{n-1}$ and 
\begin{eqnarray*}
a_{1,j}&=&\s_1 \cdots \s_{j-1} \s_j^2 \s_{j-1}^{-1} \cdots \s_1^{-1}
\end{eqnarray*}
as well as by the singular elements
\begin{eqnarray*}
X_{1,1}&=& \s_1 \t_1 \\
\t_2&=& \s_1 \s_2 \t_1 \s_2^{-1} \s_1^{-1}.
\end{eqnarray*}
$SGD_{n,1}$ is defined by the following relations:
\begin {enumerate}
\item The relations of $D_{n,1} \subset B_n$
\item The relations involving singular generators coming from the embedding
$$\phi _{n,1}: SG_{n-1} \rightarrow SGD_{n,1}:$$

\begin {eqnarray}
\s_i \t_2&=&\t_2 \s_i \mbox{ for } i \geq 2; i \neq 3 \label{Rel 1 for t1 si}\\ 
\s_3 \s_2^2 \s_3 \t_2&=& \t_2 \s_3 \s_2^2 \s_3\\
\s_3 \s_4 \s_2 \s_3 \t_2 \s_3 \s_4 \s_2 \s_3 \t_2 &=& 
\t_2 \s_3 \s_4 \s_2 \s_3 \t_2 \s_3 \s_4 \s_2 \s_3
\end{eqnarray}

\item $a_{1,i} \t_2 = \t_2 a_{1,i} \qquad$ for $i\geq 3$  \label{Rel 2 for t1 si}
\item $a_{1,2} a_{1,1} \t_2 = \t_2 a_{1,2} a_{1,1}$      
\item $\s_i X_{1,1} = X_{1,1} \s_i \qquad$ for $i\geq 3$  \label{Rel 3 for t1 si}
\item $a_{1,1} X_{1,1} = X_{1,1} a_{1,1}$  \label{Rel 4 for t1 si}
\item $X_{1,1} \s_2 a_{1,1} \s_2 = \s_2 a_{1,1} \s_2 X_{1,1}$
\item $X_{1,1} \s_2 \s_3 a_{1,2} a_{1,1} \s_3 \s_2 = \s_2 \s_3 
a_{1,2} a_{1,1} \s_3 \s_2 X_{1,1}$ \label{item_iii}
\item $X_{1,1} \s_2 \s_3 \t_2 \s_3^{-1} \s_2^{-1} = \s_2 \s_3 \t_2
\s_3^{-1} \s_2^{-1} X_{1,1}.$ \label{item_v}
\end{enumerate}
\end{theorem}

\subsection{The {\it HNN}-group-structure of the pure singular group $SGP_n$}

\begin{proposition}  \label{form of presentation}
There is a presentation for the subgroup $SGD_{n,i}$ of $SG_n$ 
in terms of the generators $a_{k,j}$ and $X_{k,j}$ as in Lemma 
\ref {coset enumeration}
so that
the relators are either of the following forms:
\begin {enumerate}
\item relators coming from the subgroup $D_{n,i} \subset B_n$.
\item The relations coming from the embedding $\Phi_{n,i}: SG_{n-i}
\longrightarrow SGD_{n,i}$.

\item \label{formiii} $X_{k,l} w= w X_{k,l}$, where $k \leq l \leq i$ and $w$ is an
element of the
pure braid group $P_n$, written in terms of the generators of $D_{n,i}$.

\item $X_{k,l} \s_j = \s_j X_{k,l}$ for some $j> i+1$ and $k \leq l \leq i$.

\item \label{formiv} $X_{k,l} \t_{i+1}= \t_{i+1} X_{k,l}, \, k \leq l \leq i$.

\item $\t_{i+1} w = w \t_{i+1}$ where $w$ is an element of the pure braid group
$P_n$.
\item
\begin{eqnarray*}
X_{l,k} \s_{i+1} \s_{i+2} \t_{i+1} \s_{i+2}^{-1} \s_{i+1}^{-1}&=&
\s_{i+1} \s_{i+2} \t_{i+1} \s_{i+2}^{-1} \s_{i+1}^{-1} X_{l,k}
\end{eqnarray*}
for some $l \leq k \leq i$.
\item \label{formlast} $X_{k,l} w X_{r,s} w^{-1} = w X_{r,s} w^{-1} X_{k,l}$
for some $k<r\leq i$, where $w$ is a word in $D_{n,i}$. 
\end{enumerate}
\end{proposition}

Clearly if we consider the subgroup $SGP_n := SGD_{n,n-1}$ relations
involving $\t_i$ and $\s_i$ no longer occur.
Hence from our proposition it immediately follows:

\begin{theorem} \label{HNN Decomposition}
Let $X$ be the collection $\{X_{i,j}, 1 \leq i \leq j \leq n-1 \}$ of the generators
of $SGP_n$ involving singularities,
and let $A$ be the collection $\{a_{i,j}, 1 \leq i \leq j \}$ of non-singular generators.

For each choice of $X_{i,j} \in X$ is $SGP_n$ isomorphic to an {\it HNN-}extension of the
subgroup $H_{i,j}$ of $SGP_n$ that is generated by all 
$x \in X - \{X_{i,j} \}$ and all $a \in A$:

$$SGP_n = \pr{ H_{i,j}, X_{i,j}\, } {\mbox{rel.} H_{i,j},\, X_{i,j} U_{i,j} =U_{i,j} X_{i,j}\,}$$ 
for some subgroups $U_{i,j}$ in $H_{i,j}$.
\end{theorem}

Hence, the group $SGP_n$ is an iterated {\it HNN-}extension of the group $P_n$. This gives us the first Betti number:

\begin{corollary} \label{first Betti number}
The first homology group with integer coefficients is:
$$
\mbox{H}_1(SGP_n,\Z) \cong \Z^{n(n-1)}.
$$
\end{corollary}

\begin{proof}
The first homology group for the pure braid group $P_n$ is well-known
to be free abelian of rank $n (n-1)/2$.
It follows e.g. immediately from the fact that 
the short exact sequence
$$
\{0\} \longrightarrow F_{n-1} \longrightarrow P_{n} 
\stackrel {\longleftarrow} {\longrightarrow} P_{n-1} \longrightarrow \{0\}
$$
splits (see e.g. \cite{Birman3}).
Here $F_{n-1}$ is the free subgroup of rank $n-1$ in $P_n$ generated by
$\{a_{1,1}, \dots, a_{1,n-1} \}$ 
(for notations confer Lemma \ref{coset enumeration}) and
the homomorphism $P_n \longrightarrow P_{n-1}$ is given by pulling out
the first strand of a pure braid. 
Therefore - by induction - the first homology group of the
group $P_n$ is free abelian of rank $(n (n-1))/2$.

If $G$ is a group which abelianization is free abelian of rank $k$ then the
abelianization of an {\it HNN-}extension has rank $k+1$. Since the
group $SGP_n$ is an iterated {\it HNN-}extension by 
Theorem \ref {HNN Decomposition} and since the cardinality of 
the set $X$ of stable letters is of size $(n(n-1))/2$ 
we get the desired result.  
\end{proof} 

We will need an additional lemma to Proposition \ref {form of 
presentation} which follows easily from geometrical considerations:

\begin{lemma} \label{cannot commute}
A relation 
$$
X_{i,j} w_1 X_{k,l} w_2 = w_1 X_{k,l} w_2 X_{i,j}
$$
with $w_1, w_2 \in P_n$ cannot occur in $SGP_n$ if $i=k$ or $i=l+1$ or
$j+1=k$ or $j=l$.
\end{lemma}

As a corollary to Theorem \ref{HNN Decomposition} and Lemma 
\ref{cannot commute} we get a presentation for the following factor group: 

\begin{corollary}
Let $N$ be the subgroup of $SGP_n$ normally generated by $P_n$.
Then $SGP_n/N$ has the presentation
$$
\pr {X \,} {X_{i,j} X_{k,l} = X_{k,l} X_{i,j} \, \mbox{ for } i \neq k, l+1, \,
\mbox{ and } j \neq l, k-1 \, }.
$$ 
\end{corollary}

\subsection{Example: The pure singular braid group on three strands}

\begin{example} \rm 
The group $SGP_3$ is generated by the elements $a_{1,1}, a_{1,2}, a_{2,2}$
as well as $X_{1,1}, X_{1,2}$ and $X_{2,2}$.

The relations are:
\begin{eqnarray*}
a_{2,2} a_{1,2} a_{2,2}^{-1} &=& a_{1,1}^{-1} a_{1,2} a_{1,1}\\
a_{2,2} a_{1,2} a_{1,1} a_{2,2}^{-1} &=& a_{1,2} a_{1,1}\\
X_{2,2} a_{2,2} &=& a_{2,2} X_{2,2}\\
X_{2,2} a_{1,2} a_{1,1} &=& a_{1,2} a_{1,1} X_{2,2}\\
X_{1,1} a_{1,1} &=& a_{1,1} X_{1,1}\\
X_{1,1} a_{1,2} a_{1,1} a_{2,2} &=& a_{1,2} a_{1,1} a_{2,2} X_{1,1}\\
X_{1,2} a_{1,1}^{-1} a_{1,2} a_{1,1}&=& a_{1,1}^{-1} a_{1,2} a_{1,1} X_{1,2}\\
X_{1,2} a_{2,2} a_{1,1} &=& a_{2,2} a_{1,1} X_{1,2}
\end{eqnarray*}
\end{example}

\section{A solution for the word problem in $SB_n$}

To give a solution to the word problem in $SB_n$ we will proceed as follows:
We know that $SB_n$ embeds into a group $SG_n$. Especially we know that
two words $w_1 s$ and $w_2 s$ in $SB_n$, where $w_1$, $w_2$ and $s$ are in
$SB_n$, are equivalent if and only if $w_1$ and $w_2$ are equivalent in
$SB_n$.

We note that - for our purposes - it is sufficient to solve the 
word problem for any two words $w_1$ and $w_2$ in $SGP_n$ with positive 
exponents for each singular generator.

Again let $X$ be the set of singular generators $\{X_{i,j}\}$ of $SGP_n$.

\begin{theorem}  \label{maintheorem}
Let $w_1 = \alpha_1 Y_1 \alpha_2 \cdots \alpha_m Y_m$ and
$w_2 = Z_1 \beta_1 \cdots Z_r \beta_r$ be two words in $SP_n \subset SGP_n$ 
with $Y_j, Z_k \in X$ and $\alpha_j, \beta _k \in P_n$.

Then $w_1 = w_2$ if and only if
the following hold:

There is a $j$ such that
\begin {eqnarray}
Z_j&=& Y_m \mbox{ and }  Z_i \neq Y_m \mbox { for } i>j 
\label{equation_1}\\
Z_{r-l} \beta_{r-l} \cdots \beta_{r-1} \beta_r Y_m \beta_r^{-1} \cdots
\beta^{-1}_{r-l} &=&
\beta_{r-l} \cdots \beta_{r-1} \beta_r Y_m \beta_r^{-1} \cdots
\beta^{-1}_{r-l} Z_{r-l} \label{equation_2}\\
& &  \mbox{ for all } r-l>j \nonumber \\
\beta_j \cdots \beta_r Y_m &=& Y_m \beta_j \cdots \beta_r \label{equation_3}\\
\alpha_1 Y_1 \alpha_2 \cdots Y_{m-1} \alpha_m&=&
Z_1 \beta_1 \cdots \beta_{j-1} \beta_j \cdots Z_r \beta_r \label{equation_4}.
\end{eqnarray}

This gives us a solution to the word problem.
\end {theorem}

\begin{proof}
By our Theorem \ref{HNN Decomposition} we know that we can regard
$SGP_n$ as {\it HNN-}extension with stable letter $Y_m$.
Thus, by Britton's lemma, if 
$w_1 = \alpha_1 Y_1 \alpha_2 \cdots \alpha_m Y_m
= w_2 = Z_1 \beta_1 \cdots Z_r \beta_r$ then there must be a $j$ satisfying 
(\ref{equation_1}) and 

\begin{eqnarray}
Y_m \beta_j Z_{j+1} \cdots Z_r \beta_r &=&
\beta_j Z_{j+1} \cdots Z_r \beta_r Y_m \label{equation_5}
\end{eqnarray}
and therefore (\ref{equation_4}).

The converse is also true.

If we now consider $Z_r$ as the stable letter with the same argument - 
Britton's Lemma and Theorem \ref{HNN Decomposition} -
we see that Equation (\ref{equation_5}) is equivalent to 
\begin{eqnarray*}
Z_r \beta_r Y_m \beta_r^{-1} &=& \beta_r Y_m \beta_r^{-1} Z_r\\
Y_m \beta_j Z_{j+1} \cdots \beta_{r-1} \beta_r &=&
\beta_j Z_{j+1} \cdots \beta_{r-1} \beta_r Y_m.
\end{eqnarray*}

Hence we will end up with Equations (\ref{equation_2}) and (\ref{equation_3}).

Since the word problem is solvable for the braid group the Equations
(\ref{equation_2}) and (\ref{equation_3}) are testable by
Lemma \ref {injective for two}.

Equation (\ref{equation_1}) is easy to test and the test for Equation
(\ref{equation_4})
is given by induction on the number of singular generators in a word.
\end{proof}

\begin{remark} It is not hard to see that along the same line one can
actually get a solution for the word problem in the whole group $SG_n$.
\end{remark}

\section{The complexity of the algorithm}

We know by the approach of Birman, Ko and Lee \cite{BKL} that the complexity
for the word problem in $B_n$ is in $O(\vert w \vert^2 n)$, where $\vert w \vert$
is the word length in terms of the generators $\s _j$ of $B_n$.

To avoid messy details and computations in the sequel we are only interested in the
complexity for the word problem for a fixed number $n$ of strands.
Our aim is to give the complexity for the word problem for the singular braid group
in terms of $\vert w \vert$ the total word length and $\vert w \vert _s$ the
number of singular generators in a word.

The pure braid group is generated by $a_{k,j}:= \s_k \s_{k-1} \cdots \s_j^2
\s_{j+1}^{-1} \cdots \s_k^{-1}$. 
Since we fixed the number of strands, the complexity for the word problem in this
group is also $O(\vert w \vert^2)$ where now $\vert w \vert$ is the word
length in terms of the new generators $a_{k,l}$.

Let $w_1$ and $w_2$ be two given words in the singular braid monoid $SB_n$.
We regard $w_1$ and $w_2$ also as elements of the group $SG_n$.

The factor group $SG_n/SGP_n$ is isomorphic to the symmetric group on $n$ elements.
Hence, to compute the right coset class of $w_1$ and $w_2$ modulo the subgroup $SGP_n$
of $SG_n$ is clearly linear in the word length of $w_1$ or $w_2$. If they are in
different classes, we are done. Otherwise we work with
$w_1 M_{i,j}^{-1}$ and $w_2 M_{i,j}^{-1}$ instead, where $M_{i,j}$ is the
representative of the right coset class of $w_1$ and $w_2$. 

Rewriting the new $w_1$ and $w_2$ in terms of the generators $X_{i,j}$ and $a_{i,j}$ of
$SGP_n$ does not change the word length or the number of singular generators. So we can assume without loss of generality
that $w_1$ and $w_2$ are already in $SGP_n$ and are given as products of the generators
$X_{i,j}$ and $a_{i,j}$.

Now let for two fixed words $w_1$ and $w_2$ in $SGP_n$ involving only positive exponents of the singular generators the values $\vert w \vert$ 
(resp. $\vert w \vert_s$) be the maximum of $\vert w_1 \vert$ and 
$\vert w_2 \vert$ (resp. the maximum of $\vert w_1 \vert_s$ and $\vert w_2 \vert_s$).

A short look at (\ref{equation_1}) - (\ref{equation_4}) gives the following:

\begin{enumerate}
\item To check (\ref{equation_1}) is linear in $O(\vert w \vert_s)$.
\item To check (\ref{equation_2}) is in $O(\vert w \vert_s \vert w \vert^2)$
since the word problem for $B_n$ is in $O(\vert w \vert^2)$.
\item To check (\ref{equation_3}) is in $O(\vert w \vert^2)$.
\end{enumerate}

By (\ref{equation_4}) we know that we have to check 
(\ref{equation_1}) - (\ref{equation_3}) at most $\vert w \vert_s$-times.

Hence we have proved:
\begin{theorem}
The complexity of the word problem in $SB_n$ with the above definitions is in
$O(\vert w \vert^2_s \vert w \vert^2)$. 
\end{theorem}

\section{The group $SG_n$ is torsion-free}

There are many proofs for the well-known fact that the braid groups $B_n$ are
torsion-free. Until most recently, however, none of them could be 
considered as being elementary.

Now Dehornoy's ordering of the braid group and especially the
interpretation of it, given in \cite{FGRRW}, yields an easy way to show this
result.

The proof of the torsion-freeness of the group $SG_n$ was announced in \cite
{FKR} but as far as we know was never proved. The proof, however, follows 
directly from our approach:

\begin{theorem}[Fenn, Keyman, Rourke]\label{torsionfree}
The group $SG_n$ is torsion-free.
\end {theorem}

\begin{proof}
First we note that by the structure theorems for {\it HNN-}groups (e.g.
\cite{MKS}, \cite{LS}) torsion must lie in
the base group. Thus, by Theorem \ref{HNN Decomposition}, torsion in the normal
subgroup $SGP_n$ of $SG_n$
must lie in $P_n$ which
is torsion free as a subgroup of $B_n$.
So $SGP_n$ is torsion free. 

The subgroup $N$ of $SG_n$ normally generated by $\tau_1 \sigma_1^{-1}$ is a
subgroup of $SGP_n$ and therefore also torsion-free.

As it easily follows from Proposition \ref{second presentation of SB} the group
$SG_n/N$ is isomorphic to $B_n$. Since $B_n$ is torsion-free a torsion element
must lie in $N$ and we are done.
\end{proof}

\subsection{The cohomological dimension of $SG_n$}

A good reference for almost all facts that we use about the
cohomological dimensions $\cd(G)$ of a group $G$ is \cite{Brown} or 
\cite{Serre}.
Especially we use Serre's Theorem that the cohomological dimension 
of a torsion-free group is equal to the one of each of the subgroups of 
finite index. Furthermore the cohomological dimension of a subgroup
must be less or equal to the cohomological dimension of the group.
Since the cohomological dimension of a free abelian group of rank $n$ is
$n$ this means that the following lemma gives us a lower bound.
The first part - for the braid group - is of course very well known.

\begin{lemma} \label{lower cd bound}
The pure braid group $P_n$ contains a free abelian subgroup of rank
$n-1$. The pure singular braid group contains a free abelian subgroup
of rank $n-1+ \lfloor n/2 \rfloor$.
\end{lemma}                      

\begin{proof}
The images of $a_{k,j}=\s_k \s_{k+1} \cdots \s_j^2 \s_{j-1}^{-1}
\cdots \s_k^{-1}, \, n-1 \geq j \geq k \geq 1,$ form a basis for
the commutator factor group of $P_n$, which is free 
abelian of rank $n (n-1)/2$.

Correspondingly the images of the $a_{i,j}$ and of
$X_{k,j}=\s_j \s_{j-1} \cdots \s_k \t_k \s_{k+1}^{-1} 
\cdots \s_j^{-1}$ generate by Lemma \ref{first Betti number}
commutator factor group of $SGP_n$, which is free abelian of rank
$n (n-1)$.

By a result of Chow \cite{Chow} the center of
the braid group $B_n$ is infinite cyclic and is generated
by $c_n:=(a_{1,n-1}\cdots a_{1,1})(a_{2,n-1} \cdots a_{2,2})
\cdots (a_{n-1,n-1})$. 
This was generalized by Fenn, Rourke and Zhu \cite{FRZ} to 
the singular braid monoid $SB_n$. For $n \geq 3$ the center is also infinite
cyclic and generated by $c_n$. 

For $n=2$ the group $B_2$ is infinite cyclic and 
$SG_2$ is free abelian of rank $2$. 
For $n=3$ the element $a_{2,2}$ and the center $c_3$ form a 
free abelian subgroup of rank $2$ in $P_3$ and  
$a_{2,2}, c_3$ and $X_{2,2}$ form a free abelian subgroup of rank 
$3$ in $SGP_3$.

Now for $n>3$ the elements $c_n$ and $a_{n-1,n-1}$ in $P_n$ 
both commute with each other and with $P_{n-2}=\{a_{k,j}, \,k,j 
\leq n-3\} \subset P_{n}$ and are independent in the commutator
factor group. 
Therefore, by induction, $P_n$ contains a free abelian subgroup
of rank $n-1$.

Correspondingly, $X_{n-1,n-1}, a_{n-1,n-1}$ and $c_n$ commute with each 
other and $SGP_{n-2}$. The claim follows for $SGP_{n}$.   
\end{proof}

Again for the braid group itself the following theorem is well-known.
(See e.g. \cite{Vassiliev} for an account to results of Arnold and Fuchs.)
We only include a proof for completeness.

\begin{theorem} \label{cohomological dimension}
The group $B_n$ has cohomological dimension $n-1$.
The group $SG_n$ has cohomological dimension $n-1+\lfloor n/2 \rfloor$.
\end{theorem}
\begin{proof}
The braid group $B_n$ is torsion-free and
since $SG_n$ is torsion-free by 
Theorem \ref{torsionfree} it is enough by Serre's Theorem to
prove the theorem for a subgroup of finite index in $SG_n$ and in $B_n$.
We choose $SGP_n$ and $P_n$ for this purpose.

First we give the argument for $B_n$. 
By Lemma \ref {lower cd bound} we already know that $n-1$ is a lower bound
for the cohomological dimension $\cd(B_n)$.

Now for a group $G$ and a normal subgroup $H$ in $G$ the relation
$\cd(G) \leq \cd(H)+\cd(G/H)$ holds.
The kernel of the natural map $P_n \longrightarrow P_{n-1}$ is
free of rank $n-1$. Therefore its cohomological dimension is $1$. 
Furthermore $P_2$ is infinite cyclic and therefore also of
cohomological dimension $1$. The result follows by induction. 

For $SGP_n$ the situation is more complicated.
The subgroup $SGP_n$ has by Theorem \ref{HNN Decomposition}
the structure of an iterated {\it HNN-}extension of
the subgroup $P_n$ in $B_n$.

We know (see \cite{Bieri}) that for an {\it HNN-}extension 
$G = \pr{H,t \,} {rel. H, t U t^{-1} = U\,}$ 
of a group $H$ with subgroup $U$ the relation
\begin{eqnarray} \label{cd of HNN}
\cd(G) &\leq& \max(\cd(H), \cd(U)+1) 
\end{eqnarray}
holds.

In the following we change the notation for reasons of simplifications.
We define $Y_{i,j}:=X_{i,j+1}$ if $i<j$ and $Y_{i,j}:=X_{j, i+1}$ if
$j<i$. This means $Y_{i,j}$ is a singular pure braid so that string
$i$ intersects string $j$ once.

For an index set $I \subset \{1, \dots, n-1 \}$ let $X_I$ be the set
of all singular generators $Y_{i,j}$ with $i, j \in I$.

We will show that the subgroup $H_I$ of $SGP_n$ generated by $P_n$ and $X_I$
has cohomological dimension less or equal to $\vert I \vert /2 + n-1$.

If $\vert I \vert =2$ then there is just one singular generator, say $Y_{j,k}$  
in $X_I$ and $H_I$ is {\it HNN-}extension of $P_n$ so by $(\ref{cd of HNN})$ 
we have $\cd (H_I) \leq \cd(P_n)+1 = n$.

If $\vert I \vert =3$ then there are three singular generators in $X_I$, say
$Y_{j,k}, Y_{k,l}$ and $Y_{j,l}$.
By Lemma \ref {cannot commute} we know that each of these three generators 
cannot commute with a word that includes one of the others.
Therefore by our structure theorem $H_I$ is an {\it HNN-}extension 
$$
\pr {SGP_n, Y_{j,k}, Y_{k,l}, Y_{j,l}\,}
    {\mbox{rel.} SGP_n, Y_{j,k} U_{j,k}=U_{j,k} Y_{j,k},
    Y_{k,l} U_{k,l}= U_{k,l} Y_{k,l},
    Y_{j,l} U_{j,l}= U_{j,l} Y_{j,l}\,}
$$
with three subgroups $U_{j,k}, U_{k,l}$ and $U_{j,l}$ of $P_n$.
Therefore all three subgroups have cohomological dimension less or equal
to $\cd(P_n)$ and thus $\cd(H_I) \leq \cd(P_n)+1 = n$.

For $I=I' \cup \{j\}$ for some $j$ we know that $Y_{j,i}, \, i \in I'$, cannot
commute by Lemma \ref {cannot commute} with a word that includes one singular
generator $Y_{k,i}$ for some $k \in I'$.
So the subgroup of $H_{I'}$ with which $Y_{i,j}$ commutes is
by our structure theorem actually a subgroup of $H_{I'-\{i\}}$.

Since by induction $\cd(H_{I'-\{i\}}) \leq n-1 + \vert I \vert/2 -1$ we know that
the {\it HNN-}extension
$$
\pr {Y_{j,i}, H_{I'}\,}{\mbox{rel. }H_{I'}, Y_{j,i} U_{j,i}=U_{j,i} Y_{j,i}} 
$$
must have cohomological dimension less or equal than 
$$max(cd(H_{I'}), cd(H_{I'-\{ i \} }+1) = n-1+ \lfloor |I|/2 \rfloor $$


Now if we successively add all other $Y_{j,k}, \, k \in I',$ to this group then 
by the same arguments we still have this upper bound for it, since 
$Y_{i,k}, \, i \in I',$ only commutes with a subgroup of $H_{I' - \{k\}}$.
\end{proof}

\begin{example} 
\rm We will show how to use the {\it HNN-}structure of 
$$SG_3 \cong \pr{B_3, \t_1\,}{\mbox{rel. }B_3, \, \t_1 \s_1= \s_1 \t_1,
\t_1 (\s_2 \s_1^2 \s_2) = (\s_2 \s_1^2 \s_2) \t_1 \,}$$ 
to compute the homology of this group. 

Let $K$ be the trefoil knot embedded in $S^3$. Let $U$ be a tubular 
neighborhood of $K$. By $C = S^3 \setminus U$ we denote the closure of the 
complement of this tubular neighborhood. Obviously $\partial C = \partial U$ 
is homeomorphic to the torus $T^2$. 

It is well known that the space $C$ is a $K(B_3,1)$-space. In fact,
the fundamental group of $C$ is isomorphic to $B_3$ and since $C$ is
the closure of the complement of the tubular neighborhood of a knot, the higher
homotopy groups are trivial (see e.g. \cite{BZ}).

Moreover, the embedding of $\partial C$ into $C$ induces an injection $i_{\sharp}$ 
from $\pi_1(\partial C,*) \cong \Z \oplus \Z$ into $B_3$. The image of 
$i_{\sharp}$ is generated by $\s_1$ and $\s_2 \s_1^2 \s_2$ 
(see figure below).

\begin{figure}[h]
\centerline{\hbox{\epsfbox{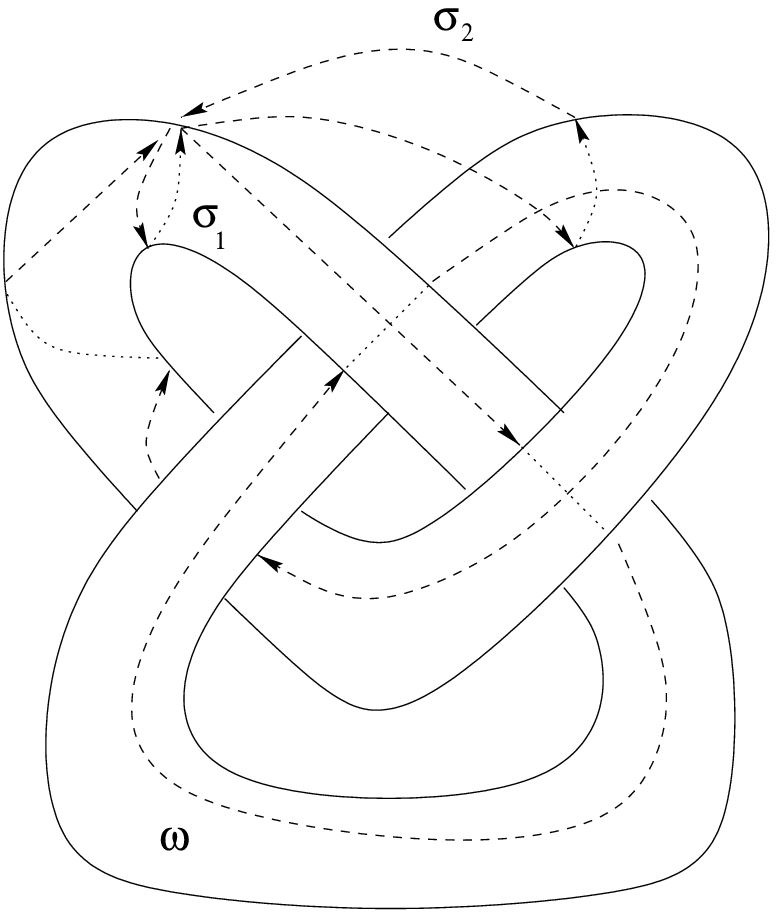}}}
\caption{$\omega = \s_2 \s_1^2 \s_2$}
\end{figure}

Thus $SG_3$ may be regarded as the {\it HNN-}extension with stable parameter 
$\t_1$ and invariant subgroup $Im(i_{\sharp})$. 

This allows us to construct a $K(SG_3,1)$-space in the following way: 
Consider the space $E = S^1 \times S^1 \times I$. Its boundary consists 
of two solid tori $T_0 = S^1 \times S^1 \times \{ 0 \}$ and 
$T_1 = S^1 \times S^1 \times \{ 1 \}$. Take a function 
$f: T \rightarrow \partial C$ which sends the longitude of $T$ to 
$\s_2\s_1^2 \s_2$ and the meridian of $T$ to $\s_1$. Attach 
both tori, $T_0$ and $T_1$, to $C$ using $f$ as an attaching map, 
in order to obtain a space 
$X = E \coprod_{f \coprod f} C$. Note that, due to the special structure 
of our attaching map, the image $\hat E$ of $E$ in $X$ is homeomorphic 
to $S^1 \times S^1 \times S^1$. 

The space $X$ is a $K(SG_3,1)$-space. Indeed, it may be easily seen 
that $\pi_1(X,*) \cong SG_3$ using the theorem of Seifert-van Kampen. 
The fact that $\pi_n(X,*) = 0$ for $n > 1$ follows from a general 
argument using covering spaces. In fact, attaching several copies of 
the universal covering space of $E$ to the universal 
covering space of $C$ in an appropriate way, yields a covering space 
of $X$ which is homotopy equivalent to a wedge 
$\bigvee_{i \in I} S^1$ of circles. Since $\bigvee_{i\in I} S^1$
 is an Eilenberg-MacLane-space, so is $X$.

Hence, we can compute the homology of $SG_3$ by calculating the homology 
of the space $X$.   

Since $X$ is three-dimensional, we immediately 
get $H_n(X) = H_n(SG_3) = 0$ for $n \geq 4$, as it follows from Theorem
\ref{cohomological dimension}.
Moreover, 
$H_0(SG_3) \cong \Z$ and $H_1(SG_3) \cong \Z \oplus \Z$ 
since the abelianization of $SG_3$ is free abelian of rank 2. 

It remains to compute $H_2(X)$ and $H_3(X)$. This will be done by using 
the Mayer-Vietoris-Sequence 

$$ \dots \rightarrow H_3(\partial C) \rightarrow H_3(C) \oplus H_3(\hat E) 
\rightarrow H_3(X) \rightarrow $$
$$ \rightarrow H_2(\partial C) \rightarrow H_2(C) \oplus H_2(\hat E) 
\rightarrow H_2(X) \rightarrow H_1(\partial C) \rightarrow \dots $$

which arises from the natural decomposition $X = C \cup \hat E$.

It is easy to see, that the map from $H_2(\partial C)$ to 
$H_2(C) \oplus H_2(\hat E)$ is injective. Therefore the map from 
$H_3(X)$ to $H_2(\partial C)$ must be the trivial map. Since 
$H_3(\partial C) = H_3(T) = 0$ and $H_3(C) = 0$ this implies that 
$H_3(X) \cong H_3(SG_3)$ is isomorphic to 
$H_3(\hat E) \cong H_3(S^1 \times S^1 \times S^1) \cong \Z$.

We are left with the case of $H_2(X)$. As in the last case, we see that the 
map from $H_1(\partial C)$ to $H_1(C) \oplus H_1(\hat E)$ is injective. 
Thus the map from $H_2(X)$ to $H_1(\partial C)$ is the trivial map. 
Since $H_2(C) = 0$ this implies, that $H_2(X)$ is obtained as a quotient 
from $H_2(\hat E)$. In fact, a close examination of the map from 
$H_2(\partial C)$ to $H_2(\hat E)$ shows that $H_2(X) \cong H_2(SG_3)$ 
must be isomorphic to $\Z \oplus \Z$.
\end{example}

\begin{remark}
Results concerning the homology of the 
infinite singular braid group $SG_{\infty}$ may be found in 
\cite{Vershinin}.
\end{remark}

\section{The singular braid monoid embeds in a group}
\label{section embedding}

In \cite{FKR} Theorem \ref{embeds in a group} is proved with the help of
geometrical considerations. For further generalizations it might be useful
to have an algebraic proof of it at hand.

On the other hand, other embedding theorems like Ores theorem -
that was for example used in in Garsides solution for the word and conjugacy
problem for the braid groups \cite{Garside} when showing that the semigroup of 
positive braids embeds into the braid groups -  do not seem to be applicable.

Using our tools, however, it is quite easy to give an algebraic proof of the
theorem:

\begin{theorem}[\cite{FKR}] 
The singular braid monoid $SB_n$ embeds into a group $SG_n$.
\end {theorem}

\begin{proof}
Let $w_1$ and $w_2$ be two different elements of $SB_n$ that have equal 
images in $SG_n$, also denoted by $w_1$ and $w_2$. 
For an element $w$ in $SB_n$ the image in the symmetric group $\Sigma_n$ under
the natural map $SB_n \longrightarrow \Sigma_n$ is the same as under the map
$SB_n \longrightarrow SG_n \longrightarrow \Sigma_n$.

Hence, by multiplying both elements with the same element in $B_n$ we can assume
that $w_1$ and $w_2$ map to the subgroup $SGP_n$ of $SG_n$. 
Therefore we can regard them as given as words
in the generators $a_{i,j}:=\s_j \s_{j+1} \cdots \s_i^2 \s_{i-1}^{-1} \cdots \s_j^{-1}$
and $X_{i,j}:=\s_j \s_{j-1} \cdots \s_i \t_i \s_{i+1}^{-1} \cdots \s_j^{-1}$,
where the $X_{i,j}$ only occur with positive exponents.

We assume that $w_1$ and $w_2$ are minimal examples with these properties, that 
means the sum of the number of singular generators in $w_1$ and $w_2$ is minimal.

We know that $B_n$ embeds into both $SB_n$ and $SG_n$ and therefore there must
be at least one singular generator in $w_1$ or $w_2$.

Now we apply Theorem \ref{maintheorem} to the two words $w_1$ and $w_2$ 
that are equal in $SGP_n$.
We have in mind that we used the fact there - which was already induced
by the embedding theorem - that for $w_1, w_2, s \in SB_n$ we
have $w_1 s = w_2 s \iff w_1 = w_2$. This cannot cause trouble here since
we assume $w_1$ and $w_2$ to be minimal examples.

Now we know by Theorem \ref{maintheorem} that $w_1$ and $w_2$ are either 
already equal in $SB_n$ or 
there are two different words $v_1$ and $v_2$ in $SB_n$ that both maps to the 
same element
in $SGP_n$ and the sum of the numbers of singular generators in $v_1$ and $v_2$
is less than the sum of the numbers of singular generators in $w_1$ and $w_2$.

Therefore the theorem follows.
\end{proof}

\section{Proofs of Theorem \ref{Pres SGD_n1} and Theorem \ref{form of presentation}}
\label{section proofs}

\begin {proof} 
{\bf of Theorem \ref{Pres SGD_n1}}

First we note that the Relations (\ref{item_iii}) and (\ref{item_v}) 
hold in $SGD_{n,1}$.


Furthermore we already have in $D_{n,1}$:
\begin {eqnarray}
a_{1,2} \s_3 a_{1,1} \s_2 &=& \s_3 a_{1,3} \s_2 a_{1,2} \label{rel_D}\\
a_{1,1} (\s_2 \s_3^2 a_{1,3} \s_2 a_{1,2})&=& (\s_2 \s_3^2 a_{1,3} \s_2
a_{1,2}) a_{1,1} \label{rel_I}\\
a_{1,3} &=& \s_3^{-1} \s_2^{-1} a_{1,1} \s_2 \s_3  \label {rel_K}.
\end{eqnarray}

By Reidemeister-Schreier we know that we get all relations in $SGD_{n,1}$ by
applying the rewriting process on all $M_j R_i M_j^{-1}$, where $M_j$ runs through
the Schreier system of right cosets of 
$SGD_{n,1}$ in $SG_n$ (see Lemma \ref {coset enumeration}) and 
$R_i$ runs through all relations of the presentation of $SG_n$. 

Since we already know a presentation for the subgroup $D_{n,1}$ of $B_n$ and since
the Schreier cosets are in $B_n$ we only have to look at the relations involving
singular generators:

\begin{enumerate}
\item Relations coming from $\s_i \t_1 = \t_1 \s_i$, for $i \neq 2$:
These are precisely the relations in (\ref{Rel 1 for t1 si}),
(\ref{Rel 2 for t1 si}), (\ref{Rel 3 for t1 si}) and 
(\ref{Rel 4 for t1 si}).


Especially from $\t_2 a_{1,3}= a_{1,3} \t_2$ and (\ref{rel_K}) it follows  that
\begin{eqnarray}
a_{1,1} (\s_2 \s_3 \t_2 \s_3^{-1} \s_2^{-1}) &=& (\s_2 \s_3 \t_2 \s_3^{-1}
\s_2^{-1}) a_{1,1} \label{rel_B}
\end{eqnarray}

\item Relations coming from $\s_2 \s_1^2 \s_2 \t_1 = \t_1 \s_2 \s_1^2 \s_2$
\begin{equation} 
\begin{array}{lrcl}
M_j, j \geq 3 \qquad & 
\s_3 \s_2^2 \s_3 \t_2 &=& \t_2 \s_3 \s_2^2 \s_3  \\ 
M_2 & \label{rel_E}
a_{1,2} a_{1,1} \t_2 &=& \t_2 a_{1,2} a_{1,1} \\
M_1 &
\s_2 ^2 a_{1,2} X_{1,1}&=& X_{1,1} \s_2 a_{1,1} \s_2\\
& \iff \s_2 a_{1,1} \s_2 X_{1,1}&=& X_{1,1} \s_2 a_{1,1} \s_2\\
M_0 &
\s_2 a_{1,1} \s_2 X_{1,1} a_{1,1}^{-1} &=& X_{1,1} a_{1,1}^{-1} \s_2^2
a_{1,2}\\
\end{array}
\end{equation}

The relation corresponding to $M_0$ follows directly from the relation corresponding to $M_1$ and (\ref{Rel 4 for t1 si}).

\item Relations coming from $\s_2 \s_3 \s_1 \s_2 \t_1 \s_2 \s_3 \s_1 \s_2 \t_1
= \t_1 \s_2 \s_3 \s_1 \s_2 \t_1 \s_2 \s_3 \s_1 \s_2$

\begin{equation}
\begin {array}{lrcl}
M_j, j \geq 4 \qquad &
\s_3 \s_4 \s_2 \s_3 \t_2 \s_3 \s_4 \s_2 \s_3 \t_2&=& \t_2 \s_3 \s_4 \s_2 \s_3 \t_2
\s_3 \s_4 \s_2 \s_3 \\
M_3 \label{obsolet1}&
\s_3 a_{1,3} \s_2 a_{1,2} X_{1,1} \s_2 \s_3 \t_2 &=& \t_2 \s_3 a_{1,3} \s_2
a_{1,2} X_{1,1} \s_2 \s_3\\
M_2 &  
\label{obsolet2}
a_{1,2} \s_3 a_{1,1} \s_2 X_{1,1} a_{1,1}^{-1} \s_2 \s_3 \t_2 &=& 
\t_2 a_{1,2} \s_3 a_{1,1} \s_2 X_{1,1} a_{1,1}^{-1} \s_2 \s_3\\
M_1&
\label{obsolet3}
\s_2 \s_3 \t_2 \s_3 a_{1,3} \s_2 a_{1,2} X_{1,1} &=& X_{1,1} \s_2 \s_3 \t_2
a_{1,2} a_{1,1} \s_3 \s_2\\
M_0 &
\label{obsolet4}
\s_2 \s_3 \t_2 a_{1,2} \s_3 a_{1,1} \s_2 X_{1,1} a_{1,1}^{-1} &=&
X_{1,1} a_{1,1}^{-1} \s_2 \s_3 \t_2 \s_3 a_{1,3} \s_2 a_{1,2}. \\
\end{array}
\end{equation}
\end{enumerate}

Now it is easy to see that the equations in (\ref{obsolet1}) - with the first
relation as an exception - follow from
(\ref{item_iii}), (\ref{rel_B}), (\ref{item_v}),
(\ref{rel_D}) and (\ref{rel_E}). 
\end{proof}

\begin{proof} 
{\bf of Theorem \ref{form of presentation}}

By Lemma \ref{coset enumeration}  
a system of Schreier right cosets of $SGD_{n,i}$ modulo $SGD_{n,i-1}$
is given by $M_{i,j} = \s_i \cdots \s_j, \quad j=i, \dots, n-1$ and
the identity. Hence,
by application of the Reidemeister-Schreier process we know that all 
relations for
$SGD_{n,i}$ can be obtained by rewriting 
$M_{i,j} R_{i-1,k} M_{i,j}^{-1}$ in terms of the generators of 
$SGD_{n,i}$, where $j=i, \dots, n-1$ and $R_{i-1,k}$ runs
through all relations of $SGD_{n,i-1}$.

More precisely: If $R_{i-1,k} = v_1 v_2 \cdots v_r$ is a relator in
$SGD_{n,i-1}$ then

\begin{equation} \label{typical_relator}
(M_{i,j} v_1 \overline{M_{i,j} v_1}^{-1})
(\overline{M_{i,j} v_1} v_2 \overline{M_{i,j} v_1 v_2}^{-1})
\cdots
(\overline{M_{i,j} v_r^{-1}} v_r M_{i,j}^{-1})
\end{equation}
is a relator in $SGD_{n,i}$ and all necessary relators in $SGD_{n,i}$ are of this
form.

By Theorem \ref{Pres SGD_n1} the claim is true for $SGD_{n,1}$.
Now assume it is true for $SGD_{n,i-1}$. We will show that all
$M_{i,j} R_{i-1,k} M_{i,j}^{-1}$ have the form that we claimed.

First we will look at the terms in (\ref{typical_relator}) coming from the
singular generators in $SGD_{n,i-1}$.
Since the $X_{k,l}, 1 \leq k \leq l \leq i-1,$ are already in the normal subgroup
$SGP_n$ of $SG_n$ it follows that $\overline {M_{i,j} X_{k,l}} = M_{i,j}$
for each right coset $M_{i,j}$.

Moreover:
\begin{eqnarray}
M_{i,j} X_{k,l} M_{i,j}^{-1} &=& X_{k,l}  \quad \mbox { for } i > l+1 \geq k+1\\
M_{i,j} X_{k, i-1} M_{i,j}^{-1} &=& X_{k,i}.
\end{eqnarray}

Especially this means: 

\begin{enumerate}
\item
A relation $X_{k,l} w = w X_{k,l}$ where $w$ is a word in
the pure braid group $P_n$ (in terms of the generators of $D_{n,i-1}$ !!!!)
yields a relation $X_{k, m} \tilde w = \tilde w X_{k,m}$ for some $m$ and
$\tilde w$ a word in the pure braid group $P_n$ (in terms of the generators
of $D_{n,i}$).

\item
A relation $X_{k,l} \s_j = \s_j X_{k,l}, j \geq i,$ leads to a relation
either of the form $X_{k,m} \s_q = \s_q X_{k,m}$ for some $m$ and
$q >i$ or
$X_{k,m} a_{i,q} = a_{i,q} X_{k,m}$ for some $m$ and $q$.

\item 
A relation of the form (\ref{formiv}) yields relations of the form (\ref{formiv})
and form (\ref{formlast}).

\item A relation $\t_i w = w \t_i$, $w \in P_n$ yields relations of either of the
following types:

We know that $\overline{M_{i,j} \t_i}=M_{i,j}$ for $j >  i$,
$\overline{M_{i,i} \t_i}= \mbox{id }$ and $\overline{\t_i}=M_{i,i}$.

Therefore we have for $j > i$:
\begin{eqnarray} 
(M_{i,j} \t_i M_{i,j}^{-1})(M_{i,j} w M_{i,j}^{-1})&=&
(M_{i,j} w M_{i,j}^{-1}) (M_{i,j} \t_i M_{i,j}^{-1})\\
\iff \qquad \qquad \qquad \t_{i+1} \tilde w&=& \tilde w \t_{i+1}
\end{eqnarray}
for some $\tilde w \in P_n$.

Since $\t_i$ commutes with $w$ so does $M_{i,i}= \s_i$. 
This means that the word $M_{i,i} w M_{i,i}^{-1}$ is equal to $w$ in $P_n$ and
thus can be transformed to it just by the relations in $P_n$.

Therefore:
\begin{eqnarray*}
(M_{i,i} \t_i) w &=& (M_{i,i} w M_{i,i}^{-1}) (M_{i,i} \t_i)\\
\iff \qquad X_{i,i} w &=& w X_{i,i}.
\end{eqnarray*}

Finally, 
\begin{eqnarray*}
(\t_i M_{i,i}^{-1})(M_{i,i} w M_{i,i}^{-1}) &=& w (\t_i M_{i,i}^{-1})\\
\iff \qquad X_{i,i} a_{i,i}^{-1} w &=& w X_{i,i} a_{i,i}^{-1}\\
\iff \qquad X_{i,i} w a_{i,i}^{-1} &=& w X_{i,i} a_{i,i}^{-1}\\
\iff \qquad X_{i,i} w&=& w X_{i,i}
\end{eqnarray*}

\item The rewriting process for the  relation
\begin{eqnarray} 
M_{i,j} X_{l,k} \s_i \s_{i+1} \t_i \s_{i+1}^{-1} \s_i^{-1} M_{i,j'}^{-1} &=&
M_{i,j} \s_i \s_{i+1} \t_i \s_{i+1}^{-1} \s_i^{-1} X_{l,k} M_{i,j'}^{-1},
\end{eqnarray}
for a suitable $j'$ and for $l \leq k \leq i-1$  
yields for $j>i+1$:
\begin {eqnarray}
X_{l,k'} \s_{i+1} \s_{i+2} \t_{i+1} \s_{i+2}^{-1} \s_{i+1}^{-1}&=&
\s_{i+1} \s_{i+2} \t_{i+1} \s_{i+2}^{-1} \s_{i+1}^{-1} X_{l,k'}
\end{eqnarray}
for some $k' \leq i$.

Furthermore we have for $j=i+1$, since $\s_{i+1} a_{i,i+1} = a_{i,i} \s_{i+1}$:
\begin{eqnarray}
X_{l,k'} \s_{i+1} a_{i,i+1} X_{i,i} \s_{i+1}^{-1} a_{i,i}^{-1} &=&
\s_{i+1} a_{i,i+1} X_{i,i} \s_{i+1}^{-1} a_{i,i}^{-1} X_{l,k'}\\
\iff X_{l,k'} a_{i,i} \s_{i+1} X_{i,i} \s_{i+1}^{-1} a_{i,i}^{-1}&=&
\label {basic_equation}
a_{i,i} \s_{i+1} X_{i,i} \s_{i+1}^{-1} a_{i,i}^{-1} X_{l,k'},
\end{eqnarray}
for some $k' \leq i$.

For $j=i$ we get:
\begin{eqnarray}
X_{l,k'} a_{i,i} \s_{i+1} X_{i,i} a_{i,i}^{-1} a_{i, i+1}^{-1} \s_{i+1}^{-1}&=&
a_{i,i} \s_{i+1} X_{i,i} a_{i,i}^{-1} a_{i,i+1}^{-1} \s_{i+1}^{-1} X_{l,k'}. 
\end{eqnarray}
This is covered by Relation (\ref{basic_equation}) and an additional relation
that we can add:
\begin{eqnarray}
X_{l,k'} a_{i,i} \s_{i+1} a_{i,i}^{-1} \s_{i+1}^{-1} a_{i,i}^{-1}&=&
a_{i,i} \s_{i+1} a_{i,i}^{-1} \s_{i+1}^{-1} a_{i,i}^{-1} X_{l,k'}.
\end{eqnarray} 
The additional relation is of the form (\ref {formiii}).

Finally for the identity as the right coset we get:
\begin{eqnarray}
X_{l,k} \t_{i+1} &=& \t_{i+1} X_{l,k}.
\end{eqnarray}
\item A relation $X_{k,l} w X_{r,s} w^{-1}= w X_{r,s} w^{-1} X_{k,l}$ becomes:
\begin {eqnarray*}
(M_{i,j} X_{k,l} M_{i,j}^{-1}) (M_{i,j} w M_{i,\tilde j}^{-1})
(M_{i, \tilde j} X_{r,s} M_{i,\tilde j}^{-1}) (M_{i, \tilde j} w^{-1}
M_{i,j}^{-1})&=&\\
(M_{i,j} w M_{i,\tilde j}^{-1})
(M_{i, \tilde j} X_{r,s} M_{i,\tilde j}^{-1}) (M_{i, \tilde j} w^{-1}
M_{i,j}^{-1}) (M_{i,j} X_{k,l} M_{i,j}^{-1})&&
\end {eqnarray*}
\begin{displaymath}
\iff X_{k, \tilde l} \tilde w X_{r, \tilde s} \tilde w ^{-1}= \tilde w 
X_{r, \tilde s} \tilde w ^{-1} X_{k,\tilde l}
\end {displaymath}
for some $\tilde j, \tilde l, \tilde s$ and a word $\tilde w$ in $D_{n,i}$.
\end{enumerate}

\end {proof}

\bibliography{../linklit}
\bibliographystyle{amsalpha}
\end{document}